\documentclass[graybox]{svmult}

\usepackage{type1cm}        % activate if the above 3 fonts are
\usepackage{makeidx}         % allows index generation
\usepackage{graphicx}        % standard LaTeX graphics tool
\usepackage{multicol}        % used for the two-column index
\usepackage{multirow}        % used for the two-column index
\usepackage[bottom]{footmisc}% places footnotes at page bottom
\usepackage{adjustbox}       % flexible box to scale tables

\usepackage{newtxtext}       % 
\usepackage[varvw]{newtxmath}       % selects Times Roman as basic font
\newcommand{\STAB}[1]{\begin{tabular}{@{}c@{}}#1\end{tabular}}

\usepackage{ulem}
\definecolor{darkgreen}{RGB}{0,120,0}
\makeindex             % used for the subject index
\usepackage{xcolor}

\begin{document}

\title*{Towards a parallel Schwarz solver framework  for virtual elements using GDSW coarse spaces}
\titlerunning{Parallel GDSW solvers for virtual elements} %for an abbreviated version of
\author{Tommaso Bevilacqua\orcidID{0000-0002-0828-5018} and\\ Axel Klawonn\orcidID{0000-0003-4765-7387} and\\ Martin Lanser\orcidID{0000-0002-4232-9395} and\\ Adam Wasiak\orcidID{0000-0002-5705-9754}}
\institute{Tommaso Bevilacqua (Division of Mathematics), Axel Klawonn (Division of Mathematics and Center for Data and Simulation Science), Martin Lanser (Division of Mathematics and Center for Data and Simulation Science), Adam Wasiak (Division of Mathematics) \at University of Cologne, Cologne, \newline
\email{[tommaso.bevilacqua, axel.klawonn, martin.lanser, adam.wasiak]@uni-koeln.de}}
\maketitle

\abstract{The Virtual Element Method (VEM) is used to perform the discretization of the Poisson problem on polygonal and polyhedral meshes. This results in a symmetric positive definite linear system, which is solved iteratively using overlapping Schwarz domain decomposition preconditioners, where to ensure robustness and parallel scalability a second level has to be employed. The construction and numerical study of two-level overlapping Schwarz preconditioners with variants of the GDSW (Generalized Dryja-Smith-Widlund) coarse space are presented here. Our PETSc-based parallel implementation of GDSW and variants, combined with the Vem++ library, represent the first parallel application of these GDSW preconditioners to VEM. Numerical experiments in 2D and 3D demonstrate scalability of our preconditioners up to 1\,000 parallel cores for VEM discretizations of degrees $k=1,2$.}

\section{Introduction}
\label{intro}

Overlapping Schwarz methods are one of the most extensively studied and widely used domain decomposition approaches, both in theory and practice. A key component of any two-level overlapping Schwarz preconditioner is the coarse space, which is essential for ensuring robust convergence with respect to certain parameters as well as numerical and parallel scalability as the number of subdomains increases.

The Generalized Dryja-Smith-Widlund (GDSW) coarse space, introduced in \cite{dohrmann_family_2008,dohrmann_less_regular_2008}, has become a standard and robust choice for domain decomposition preconditioners due to its ease of implementation, wide applicability, and good convergence properties across diverse problem classes and discretizations. Its construction is based on decomposing the interface between subdomains into vertices, edges, and faces, defining coarse basis functions for each of them, and finally exploiting problem dependent extension operators to the interior of subdomains. 
Despite its robustness, the dimension of this coarse space is known to increase rapidly when scaled to many subdomains, especially for 3D problems.

Several strategies have been introduced to reduce the size of this coarse space, for example, GDSW* \cite{heinlein2025monolithicblockoverlappingschwarz,hochmuth_parallel_2020} and Reduced GDSW (RGDSW) \cite{dohrmann_design_2017}. Rather than using separate basis functions for each interface component, in particular, GDSW* combines the edge contributions with the vertex contributions while keeping the face contributions separate. In contrast, RGDSW adds both face and edge contributions to the contributions of the vertices connected to them. This can lead to significant reductions in coarse problem size, especially in three dimensions where the number of interface components grows rapidly.

While various two-level Overlapping Schwarz methods have been studied for the Finite Element Method (FEM), its application to the Virtual Element Method (VEM) is more recent \cite{calvo_overlapping_2019, calvo_virtual_2020, calvo_overlapping_2024}. VEM is a FEM developed for polygonal and polyhedral meshes, which is increasingly attractive for practical applications due to the flexibility in mesh generation and the potential for non-conforming discretizations \cite{antonietti2022virtual,mengolini2019engineering,salama2025adaptive}. 

The goal of this work is to conduct a numerical study of two-level additive Schwarz preconditioners employing the GDSW, GDSW*, and RGDSW coarse space variants for VEM discretizations of the Poisson problem. We investigate the scalability of both approaches on randomly seeded Voronoi meshes in 2D and 3D, and compare their performance for low order polynomial degrees ($k=1,2$) and parallel core counts (up to 1\,000 cores). Our implementation combines the Vem++ library \cite{dassi_library_2025}, for constructing the VEM matrices, with our own  parallel PETSc-based code of the different GDSW variants which is part of our in-house FE2TI code; see also \cite{bevilacqua_highly_2025} for a different application. Previous work that combined the Vem++ library with other block-preconditioners and the Balancing Domain Decomposition by Constraints (BDDC) method can be found in \cite{barnafi2023parallel,bevilacqua2024bddc,dassi2020parallel,dassi2022robust}. Our results demonstrate the expected scalability of the different variants of GDSW. Furthermore, GDSW* and RGDSW achieve significant reductions in coarse problem size while maintaining competitive convergence rates.

\section{Virtual Element Method}

We consider the weak formulation of the Poisson problem with inhomogeneous Dirichlet boundary conditions in a bounded polygonal domain $\Omega \subset \mathbb{R}^d$, $d=2,3$:
Find $w \in H^1_0(\Omega)$ such that
\begin{equation*}
a(w,v) = F(v) \quad \forall v \in H^1_0(\Omega),
\end{equation*}
where $a(w,v) = \int_\Omega \nabla w \cdot \nabla v \, dx$ and $F(v) = \int_\Omega f v \, dx - \int_\Omega \nabla \tilde{u} \cdot \nabla v \, dx$, with $\tilde{u} \in H^1(\Omega)$ satisfying $-\Delta \tilde{u} = 0$ in $\Omega$ and $\tilde{u}|_{\partial\Omega} = g$, and $u = \tilde{u} + w$, where $f \in L^2(\Omega)$ is the source term and $g \in H^{1/2}(\partial\Omega)$ is the given boundary data.
The Virtual Element Method (VEM)  extends finite element techniques to general polygonal/polyhedral meshes. Given a polygonal/polyhedral tessellation $\mathcal{T}_h$ with a maximum element size $h$ of $\Omega$, on each element $E \in \mathcal{T}_h$, the local space $V_h^E$ of order $k \geq 1$ is defined implicitly through polynomial projections, with degrees of freedom typically consisting of vertex values and edge moments.

The bilinear form $a_h(\cdot,\cdot)$ is constructed as $a_h(v_h,w_h) = a(\Pi^\nabla_k v_h, \Pi^\nabla_k w_h) + s_h(v_h - \Pi^\nabla_k v_h, w_h - \Pi^\nabla_k w_h) $, where $\Pi^\nabla_k$ is the $H^1$-projection onto polynomials of degree$\,\leq\! k$, $a(\cdot,\cdot)$ is the exact bilinear form, and $s_h(\cdot,\cdot)$ is a stability term that is computable from the degrees of freedom. We use a variant of the d-recipe stabilization for the VEM discretization in both 2D \cite{mascotto_ill_conditioning_2018} and 3D \cite{dassi_high_order_2018}. The right hand side term is given by $F_h(v_h) = F(\Pi^0_k v_h)$ where $\Pi^0_k$ is the $L^2$-projection onto polynomials of degree $\leq k$, which is computable from the degrees of freedom as well.

The VEM approximation $u_h \in V_h \subset H^1_0(\Omega)$ satisfies
\begin{equation*}
a_h(u_h,v_h) = F_h(v_h) \quad \forall v_h \in V_h,
\end{equation*}
where $a_h(\cdot,\cdot)$ and $F_h(\cdot)$ are consistent and stable approximations of $a(\cdot,\cdot)$ and $F(\cdot)$. For further details on the necessary mesh assumptions and the construction of the VEM spaces, we refer to \cite{ahmad_equivalent_2013,beirao_basic_2013}. In the following, we will assume the polynomial degree of the Virtual Element spaces to be $k=1$ or $k=2$.

\section{Overlapping Schwarz Method with GDSW}

The additive two-level GDSW preconditioner for the VEM discretization is given by 
\[
M_\text{GDSW}^{-1} = \Phi K_0^{-1} \Phi^T + \sum_{i=1}^N R_i^T K_i^{-1} R_i,
\]
where $R_i$ are restrictions to overlapping subdomains $\Omega_i$ and the local stiffness matrices are extracted from $K$ by $K_i = R_i K R_i^T$, for $i = 1, . . . , N,$. Additionally, $R_0:=\Phi^T$ is the restriction operator from the VE space to the coarse space, and $K_0 = R_0 K R_0^T$ is the coarse matrix.

The matrix $\Phi$ contains the GDSW coarse space basis functions $\phi^i$ for $i = 1,...,N_\phi$ in its columns, where $N_\phi$ is their total number. These basis functions are constructed by decomposing the interface degrees of freedom into vertices $\mathcal{V}$, edges $\mathcal{E}$, and faces $\mathcal{F}$, and defining for each component $\mathcal{P} \in \{\mathcal{V}, \mathcal{E}, \mathcal{F}\}$ its corresponding interface function $\phi_\mathcal{P}$ to be 1 on all nodes belonging to $\mathcal{P}$ and 0 on all the other interface components. As a result, the sum of all coarse basis functions build a partition of unity on the interface. The interior part of the coarse function is then given by $\phi_I = -K_{II}^{-1} K_{I\mathcal{P}} \phi_\mathcal{P}$. If $\mathcal{P}$ is the $i$-th component after enumerating the vertices, edges, and faces, then $\phi^i = (\phi_I, \phi_\mathcal{P})$. Since $K$ is symmetric positive definite in our model problem, the coarse basis functions build a partition of unity on the whole domain $\Omega$ now.\\

%\subsection{Reduced GDSW (RGDSW)}

\textbf{Reduced GDSW (RGDSW).} The RGDSW coarse space is a variant of the standard GDSW coarse space, designed to reduce the dimension of the coarse problem while maintaining robust convergence properties.
For each vertex $\mathbf{v}$ of the domain decomposition, we define the corresponding RGDSW basis function $\phi^{\mathbf{v}}$ as follows. Let $\mathcal{F}_{\mathbf{v}}$ denote the (open) set of all face degrees of freedom connected to $\mathbf{v}$ across different subdomains, $n_{\mathcal{F}}$ the number of subdomains sharing a face $\mathcal{F}$, and $\mathcal{E}_{\mathbf{v}}$ denote the (open) set of all edge degrees of freedom connected to vertex $\mathbf{v}$.
Lastly, let $\mathcal{P}_{\mathbf{v}}:= \{\mathbf{v}\}\cup \mathcal{E}_{\mathbf{v}} \cup \mathcal{F}_{\mathbf{v}}$ be the union of these sets.
The basis function $\phi^{\mathbf{v}}$ is constructed by prescribing the following values on the interface degrees of freedom:
\begin{equation*}
\phi_{\mathcal{P}_{\mathbf{v}}}(x) = \begin{cases}
1 & x=\mathbf{v}, \\
0.5 & x \in  \mathcal{E}_{\mathbf{v}}, \\
\tfrac{1}{n_{\mathcal{F}_\mathbf{v}}} & x \in  \mathcal{F}_{\mathbf{v}}, \\
0 & \text{else},
\end{cases}
\end{equation*}

Once these interface values are specified, the interior part of the basis function is computed by solving a local problem, analogous to the GDSW case.\\

%\subsection{GDSW*}

\textbf{GDSW*.} For 3D problems the GDSW* coarse space is an intermediate variant between GDSW and RGDSW, designed to reduce the dimension of the coarse problem by combining only vertex and edge contributions, while treating faces separately as in the standard GDSW approach. Since the interface decomposition of 2D problems consists only of vertices and edges (and no faces) the GDSW* and RGDSW coarse spaces coincide in this case.

The GDSW* coarse space consists of two types of basis functions:
\begin{itemize}
\item Vertex-edge basis functions: For each vertex $\mathbf{v}$ of the domain decomposition, we define a basis function $\phi^{\mathbf{v}}$ that combines the vertex with its connected edges. Let $\mathcal{E}_{\mathbf{v}}$ denote the (open) set of all edge degrees of freedom connected to vertex $\mathbf{v}$, and let $\bar{\mathcal{P}}_{\mathbf{v}}:= \{\mathbf{v}\}\cup \mathcal{E}_{\mathbf{v}}$ be their union. The basis function is constructed by prescribing:
\begin{equation*}
\phi_{\bar{\mathcal{P}}_{\mathbf{v}}}(x) = \begin{cases}
1 & x=\mathbf{v}, \\
0.5 & x \in  \mathcal{E}_{\mathbf{v}}, \\
0 & \text{else}.
\end{cases}
\end{equation*}
\item Face basis functions: For each face $\mathcal{F}$ of the domain decomposition interface, we define a separate basis function $\phi_{\mathcal{F}}$ exactly as in the standard GDSW approach. The function is set to 1 on all degrees of freedom belonging to $\mathcal{F}$ and 0 on all other interface components.
\end{itemize}
The interior part of each basis function is again computed analogously to the GDSW case. This construction yields a coarse space that is larger than RGDSW but smaller than the full GDSW coarse space, as it eliminates the separate edge basis functions while retaining the face contributions. Let us note that in all three cases (GDSW, RGDSW, GDSW*) a partition of unity is built on $\Omega$ and only the number of coarse basis functions varies.

\section{Numerical Results}

\begin{figure}[b]
\centering
\begin{minipage}{0.48\textwidth}
    \centering
    \includegraphics[width=\textwidth,trim={23cm 20cm 23cm 20cm},clip]{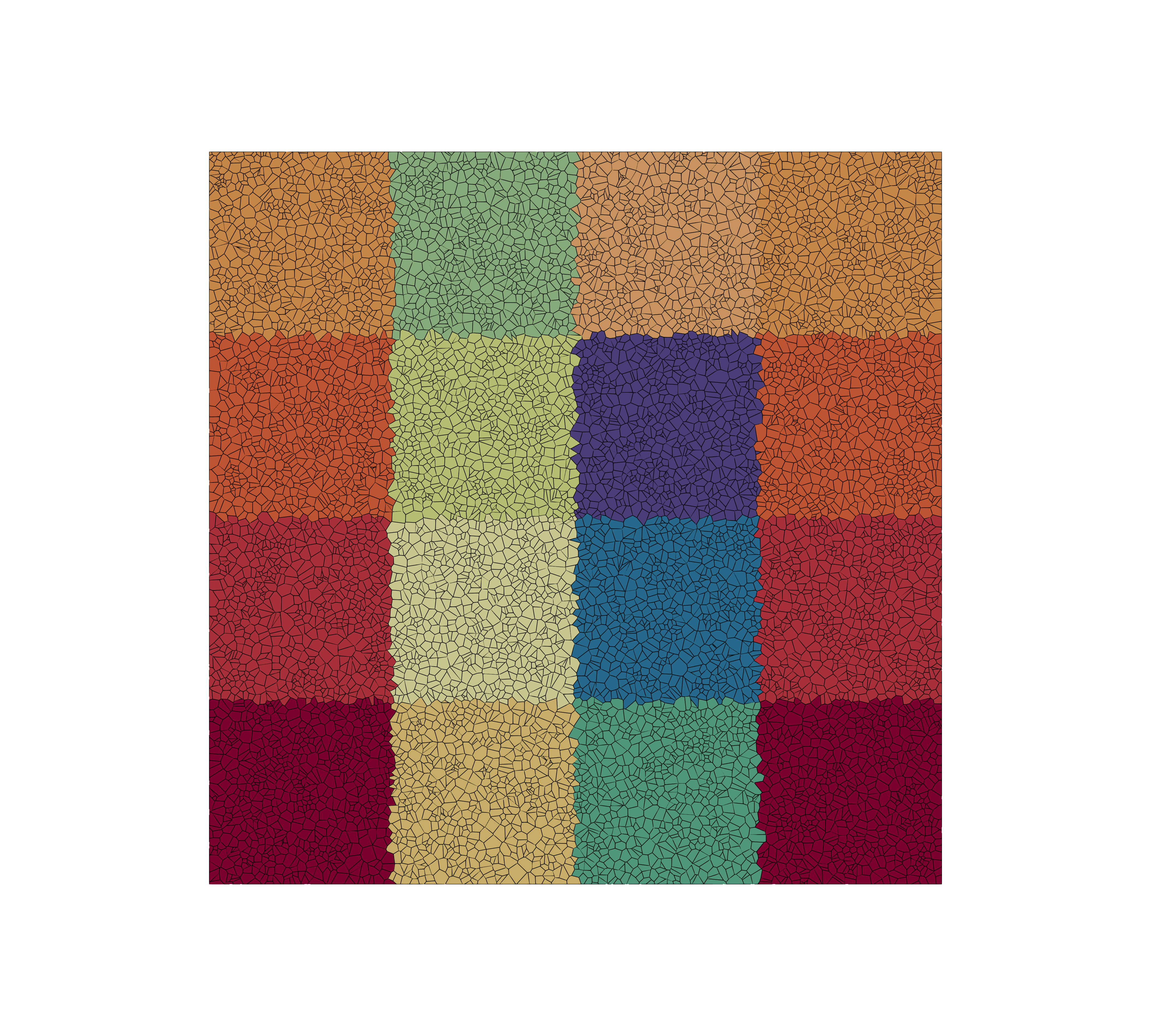}
\end{minipage}
\hfill
\begin{minipage}{0.48\textwidth}
    \centering
    \includegraphics[width=\textwidth,trim={12.2cm 2cm 15.5cm 8.5cm},clip]{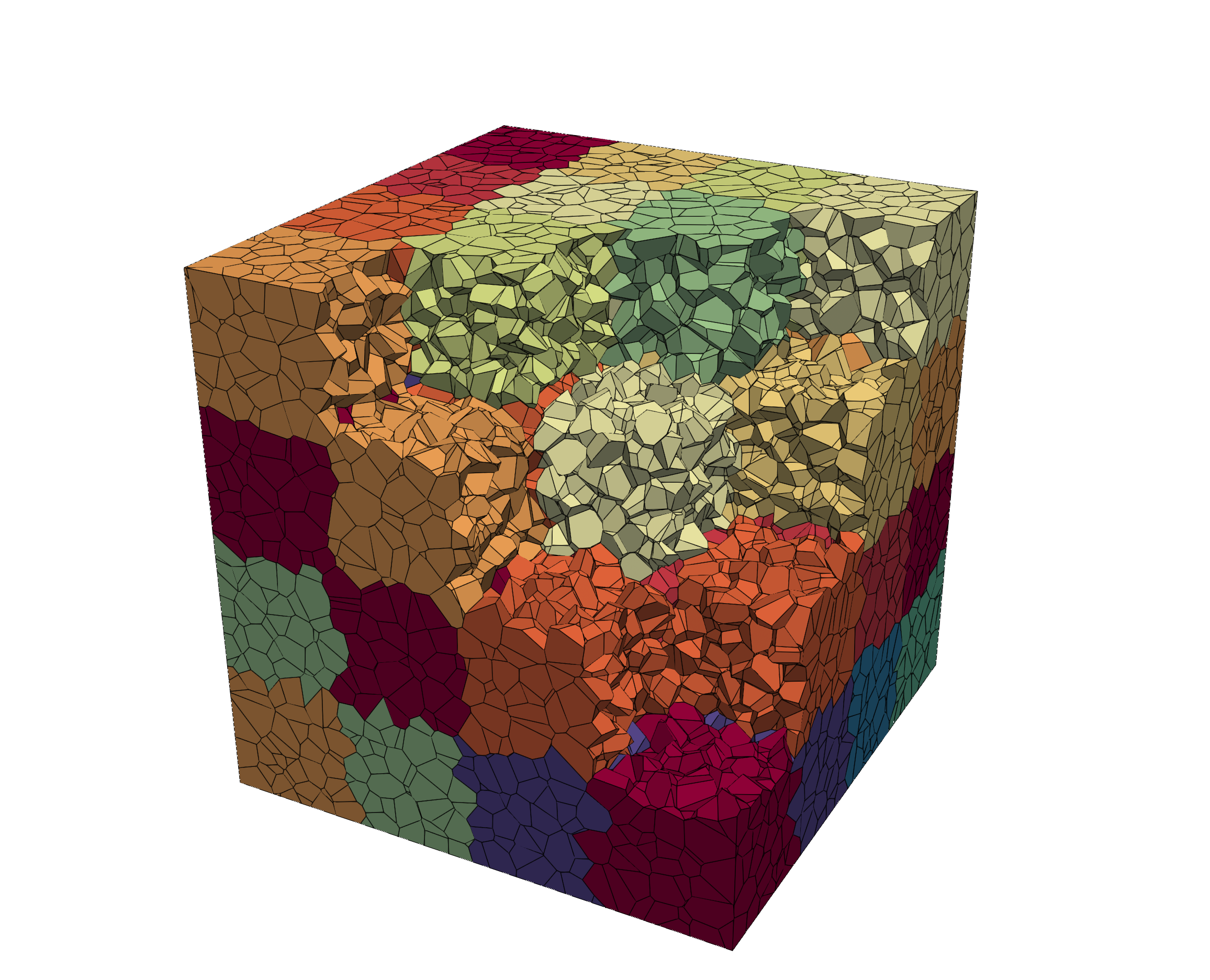}
\end{minipage}
\caption{Domain decomposition of randomly seeded Voronoi meshes: (left) 2D decomposition into 4$\times$4 subdomains of a Voronoi mesh with 10\,000 elements; (right) 3D decomposition into 4$\times$4$\times$4 subdomains of a Voronoi mesh with 16\,000 elements.}
\label{fig:due_figure}
\end{figure}

We present numerical results for the solution of the Poisson problem on the unit square in 2D and on the unit cube in 3D using the Virtual Element Method (VEM) of degree $k=1,2$. The computational domains are discretized using randomly seeded Voronoi meshes, as illustrated in Figure~\ref{fig:due_figure}, which shows the domain decomposition of a 2D mesh with 10\,000 elements and a 3D mesh with 16\,000 elements. In 2D, the source term is $f(x,y)= -6x - 2$ and the boundary function $g(x,y) = x^3+y^2$. In 3D, the source term is $f(x,y,z) = -6x - 2 - 18z$ and the boundary function is $g(x,y,z) = x^3 + y^2 + 3z^3$.

For the numerical experiments we use the Vem++ Library \cite{dassi_library_2025}, an open-source C++ library designed for solving Partial Differential Equations (PDEs) using VEM, in which we use the d-recipe stabilization found under the option \texttt{STIFFC1RECIPESTAB}. The linear system is solved by GMRES with a relative tolerance of $10^{-8}$ (of the unpreconditioned error) as stopping criterion, employing as preconditioner a two-level restricted additive Schwarz method with an overlap of 1 with the GDSW coarse space. We used the GDSW implementation of our inhouse PETSc-based FE2TI library, that allows us to compare the classical coarse space construction with the GDSW* and the RGDSW variants. 

%\subsection{Strong Scaling}

\textbf{Strong Scaling.} Table~\ref{tab:results} presents the strong scaling analysis for both 2D and 3D problems on Voronoi meshes. The results demonstrate excellent scalability of the GDSW approach across all tested configurations. For the 2D case (10\,000 elements), iteration counts remain stable between 21--27 for both $k=1$ and $k=2$ as the number of cores increases from 16 to 400. Similarly, for the 3D case (16\,000 elements), GDSW maintains stable iterations with core counts ranging from 64 to 1\,000.

The RGDSW coarse space achieves approximately 50\% reduction in coarse problem dimension compared to GDSW while maintaining competitive iteration counts. This reduction becomes increasingly important for large-scale parallel computations, demonstrating the effectiveness of the reduced coarse space construction for both 2D and 3D problems.
In 3D, the GDSW* reduces the coarse space size by approximately 25\% achieving the same iteration counts as GDSW.

\begin{table}[t]
\centering
\begin{minipage}{0.41\textwidth}
\centering
\begin{adjustbox}{width=\textwidth}
\begin{tabular}{|c|r||r|r||r|r|}
\cline{2-6}
\multicolumn{1}{c|}{} & nDoFs & \multicolumn{2}{c||}{20\,002} & \multicolumn{2}{c|}{60\,003} \\
\cline{2-6}
\multicolumn{1}{c|}{} & Degree & \multicolumn{2}{c||}{$k=1$} & \multicolumn{2}{c|}{$k=2$} \\
\cline{2-6}
\multicolumn{1}{c|}{} & nCores & Coarse & it & Coarse & it \\
\cline{2-6} 
\noalign{\vskip 0.6ex}
\cline{1-6}
\multirow{4}{*}{\STAB{\rotatebox[origin=c]{90}{GDSW}}}
& 16 & 45 & 22 & 51 & 27 \\
& 64 & 226 & 22 & 259 & 25 \\
& 256 & 1\,005 & 20 & 1\,155 & 24 \\
& 400 & 1\,587 & 19 & 1\,843 & 23 \\
\hline
\hline
\multirow{4}{*}{\STAB{\rotatebox[origin=c]{90}{RGDSW}}}
& 16 & 21 & 24 & 24 & 29 \\
& 64 & 112 & 24 & 131 & 29 \\
& 256 & 515 & 22 & 600 & 26 \\
& 400 & 820 & 21 & 978 & 24 \\
\hline
\end{tabular}
\end{adjustbox}
\vspace{0.5ex}
\small 2D Voronoi mesh (10\,000 elements)
\end{minipage}\hspace{6mm}
\begin{minipage}{0.41\textwidth}
\centering
\begin{adjustbox}{width=\textwidth}
\begin{tabular}{|c|r||r|r||r|r|}
\cline{2-6}
\multicolumn{1}{c|}{} & nDoFs & \multicolumn{2}{c||}{104\,582} & \multicolumn{2}{c|}{450\,321} \\
\cline{2-6}
\multicolumn{1}{c|}{} & Degree & \multicolumn{2}{c||}{$k=1$} & \multicolumn{2}{c|}{$k=2$} \\
\cline{2-6}
\multicolumn{1}{c|}{} & nCores & Coarse & it & Coarse & it \\
\cline{2-6} 
\noalign{\vskip 0.6ex}
\cline{1-6}
\multirow{4}{*}{\STAB{\rotatebox[origin=c]{90}{GDSW}}}
& 64 & 1\,018 & 18 & 1\,257 & 24 \\
& 216 & 4\,048 & 18 & 5\,265 & 24 \\
& 512 & 9\,949 & 18 & 13\,515 & 25 \\
& 1\,000 & 19\,094 & 17 & 27\,080 & 24 \\
\hline
\hline
\multirow{4}{*}{\STAB{\rotatebox[origin=c]{90}{GDSW*}}}
& 64 & 761& 18 & 953 & 24 \\
& 216 & 2\,827 & 19 & 3\,760 & 24 \\
& 512 & 6\,823 & 18 & 9\,470 & 25 \\
& 1\,000 & 13\,350 & 17& 18\,954 & 24 \\
\hline
\hline
\multirow{4}{*}{\STAB{\rotatebox[origin=c]{90}{RGDSW}}}
& 64 & 527 & 19 & 591 & 24 \\
& 216 & 2\,072 & 20 & 2\,366 & 25 \\
& 512 & 5\,117 & 22 & 6\,074 & 26 \\
& 1\,000 & 10\,218 & 25 & 12\,578 & 28 \\
\hline
\end{tabular}
\end{adjustbox}
\vspace{0.5ex}
\small 3D Voronoi mesh (16\,000 elements)
\end{minipage}
\caption{Strong scaling analysis of the VEM solver on randomly seeded Voronoi meshes. Left: 2D Voronoi mesh with 10\,000 elements; Right: 3D Voronoi mesh with 16\,000 elements. Each subtable reports, for different numbers of cores and polynomial degrees, the coarse grid size (Coarse) and the number of iterations (it) for GDSW, GDSW* (3D only), and RGDSW coarse spaces.}
\label{tab:results}
\end{table}

\textbf{Weak Scaling.} Table~\ref{tab:weak_scaling} presents a weak scaling analysis in 3D, where coarse space growth is more significant. In order to maintain a comparable load across the subdomains, for these tests, we used a structured hexahedral mesh keeping the number of elements for each subdomain fixed to 125. We can observe that the iteration counts remain stable across all methods as cores increase from 64 to 1\,000 and for both the degrees $k=1,2$. Both GDSW* and RGDSW maintain similar iteration counts to GDSW while reducing the coarse space size, demonstrating their effectiveness for large-scale parallel computations.

\begin{table}[b]
\centering
\begin{adjustbox}{width=0.70\textwidth}
\begin{tabular}{|c|r|r||r|r||r|r||r|r|}
\cline{4-9}
\multicolumn{1}{c}{} & \multicolumn{2}{c|}{} & \multicolumn{2}{c||}{GDSW} & \multicolumn{2}{c||}{GDSW*} & \multicolumn{2}{c|}{RGDSW} \\
\cline{2-9}
\multicolumn{1}{c|}{} & nCores & Dofs & Coarse & it & Coarse & it & Coarse & it \\
\cline{2-9}
\noalign{\vskip 0.6ex}
\cline{1-9}
\multirow{4}{*}{\STAB{\rotatebox[origin=c]{90}{$k=1$}}}
& 64 & 9\,261 & 279 & 16 & 171 & 16 & 27 & 15 \\
& 216 & 29\,791 & 1\,115 & 18 & 665 & 18 & 125 & 18 \\
& 512 & 68\,921 & 2\,863 & 19 & 1\,687 & 19 & 343 & 20 \\
& 1\,000 & 132\,651 & 5\,859 & 19 & 3\,429 & 20 & 729 & 21 \\
\hline
\hline
\multirow{4}{*}{\STAB{\rotatebox[origin=c]{90}{$k=2$}}}
& 64 & 68\,921 & 279 & 20 & 171 & 20 & 27 & 20 \\
& 216 & 226\,981 & 1\,115 & 23 & 665 & 23 & 125 & 24 \\
& 512 & 531\,441 & 2\,863 & 25 & 1\,687 & 25 & 343 & 27 \\
& 1\,000 & 1\,030\,301 & 5\,859 & 27 & 3\,429 & 27 & 729 & 29 \\
\hline
\end{tabular}
\end{adjustbox}
\caption{Weak scaling analysis of the VEM solver on 3D hexahedral meshes. The table reports, for different numbers of cores and polynomial degrees, total number of degrees of freedom (Dofs), the dimension of the coarse space (Coarse) and the number of iterations (it) for GDSW, GDSW*, and RGDSW coarse spaces.}
\label{tab:weak_scaling}
\end{table}

\section{Conclusions}

In this work, we have presented numerical scalability studies of two-level restricted additive Schwarz preconditioners for the Virtual Element Method applied to the Poisson problem in both 2D and 3D. 
Future work will focus on conducting a more detailed analysis including timings and computational efficiency. We plan to develop the theoretical framework for the convergence analysis of these preconditioners applied to VEM discretizations of the Poisson problem. We note that the above described construction of the GDSW coarse space types is only valid for VEM discretizations of polynomial orders $k=1,2$ in three dimensions. This is due to the construction of the virtual degrees of freedom. % in the interior of the elements. 
Setting the values of the basis functions to constant values for all face degrees of freedom does not accurately reflect the null space of the coarse problem. We already developed an adjusted construction for all three coarse space types for higher polynomial orders which will be presented in an upcoming publication. Furthermore, we intend to extend this study to other model problems of practical interest, particularly linear elasticity, to evaluate the robustness and effectiveness of these domain decomposition techniques across a wider class of applications.

\vspace{1em}

\begin{acknowledgement}
This project has received funding from the Deutsche Forschungsgemeinschaft (DFG) as part of the Forschungsgruppe (Research Unit) 5134 ``Solidification Cracks During Laser Beam Welding -- High Performance Computing for High Performance Processing'' under DFG project number 434946896. The authors gratefully acknowledge the scientific support and HPC resources provided by the Erlangen National High Performance Computing Center (NHR@FAU) of the Friedrich-Alexander-Universit\"at Erlangen-N\"urnberg (FAU) under the NHR project k109be10. NHR funding is provided by federal and Bavarian state authorities. NHR@FAU hardware is partially funded by the German Research Foundation (DFG) - 440719683. We would like to thank Franco Dassi from the University of Milano-Bicocca and the other developers of the Vem++ library. %{\TB [https://sites.google.com/unimib.it/vemwebsite/home-page?authuser=0]}.}
\end{acknowledgement}
\ethics{Declaration of AI usage}{In the course of preparing this work, the authors Tommaso Bevilacqua and Adam Wasiak used the current models from OpenAI and Anthropic (GPT-5, Claude 4.5 Sonnet and Haiku) to assist with typesetting the manuscript. Following the use of these tools, all authors thoroughly reviewed and revised the content as needed, taking full responsibility for the final version of the publication.}

\vspace{-2mm}

% BibTeX users please use
\bibliographystyle{spmpsci}
\bibliography{refer.bib}

\end{document}